\documentclass[12pt]{article}
\usepackage{theorem}
\newtheorem{lemma}{Lemma}
\newtheorem{proposition}[lemma]{Proposition}
\newtheorem{theorem}[lemma]{Theorem}
\newtheorem{corollary}[lemma]{Corollary}
{\theorembodyfont{\upshape}}
{\theorembodyfont{\upshape}\newtheorem{remark}[lemma]{Remark}}
{\theorembodyfont{\upshape}\newtheorem{example}[lemma]{Example}}
{\theorembodyfont{\upshape}}
{\theorembodyfont{\upshape}}
\newcommand{\R}{{\bf R}}

\newcommand{\rme}{{\rm e}}
\newcommand{\rmd}{\,{\rm d}}
\newcommand{\sig}{\sigma}
\newcommand{\alp}{\alpha}

\newcommand{\gam}{\gamma}
\newcommand{\lam}{\lambda}

\newcommand{\eps}{\varepsilon}

\newcommand{\Ome}{\Omega}
\newcommand{\lap}{{\Delta}}
\newcommand{\gra}{\nabla}

\newcommand{\norm}{\Vert}

\newcommand{\Proof}{\underbar{Proof}{\hskip 0.1in}}

\newcommand{\lin}{{\rm lin}}
\newcommand{\Schrodinger}{Schr\"odinger }
\title{SPECTRAL STABILITY\\ OF THE NEUMANN LAPLACIAN}
\author{V.~I.~Burenkov,\footnote{Supported by the INTAS grant 99-01080 and by the Russian Foundation
for Basic Research grant 99-01-00868.}\,~~ E.~B.~Davies}
\date{May 2001}
\begin{document}
\maketitle
%%%%%%%%%%%%%%%%%%%%%
\begin{abstract}
We prove the equivalence of Hardy- and Sobolev-type inequalities, certain
uniform bounds on the heat kernel and some spectral regularity
properties of the Neumann Laplacian associated with an arbitrary
region of finite measure in Euclidean space. We also prove that if
one perturbs the boundary of the region within a uniform H\"older
category then the eigenvalues of the Neumann Laplacian change by a
small and explicitly estimated amount.
\par
keywords: Neumann Laplacian, Sobolev inequalities, Hardy inequalities,
spectral stability, H\"older continuity.
\end{abstract}
%%%%%%%%%%%%%%%%%%%%%%%%
\section{Introduction}
\par
Let $\Omega$ be an arbitrary region in $\R^N$ and let us define the
Neumann Laplacian to be the non-negative self-adjoint operator
$H=-\lap_N$ acting in $L^2(\Ome)$ and associated with the quadratic
form
\begin{equation}
Q(f)=\left\{ \begin{array}{ll}
\int_\Ome |\gra f|^2 \rmd^N x & \mbox{if $f\in W^{1,2}(\Ome)$}\\
+\infty & \mbox{otherwise}
\end{array}\right.
\label{B}
\end{equation}
as described in \cite[Section 4.4]{STDO}. It is well known that if
$\Ome$ is bounded with continuous boundary $\partial
\Ome$ then $H$ has compact resolvent since the embedding
$W^{1,2}(\Omega)\subset L^2(\Omega)$ is compact, \cite{Bur}.
However, in general the spectrum of $H$ may be quite wild, even for
bounded regions in $\R^2$, \cite{HS}. These phenomena are not well
understood, with the result that the Neumann Laplacian is far less
studied than the Dirichlet Laplacian.

In this paper we prove a number of general results concerning the
spectral behaviour of the Neumann Laplacian. We start by
investigating the relationship between Hardy-type and Sobolev-type inequalities
for arbitrary regions of finite inradius. We then
establish the equivalence of Sobolev-type inequalities to some spectral
properties of the Neumann Laplacian. The results apply in particular to bounded
regions with H\"older continuous boundaries.

Even if one knows that the spectrum is discrete, the numerical
computation of the eigenvalues by the finite element or other
methods depends upon the assumption that if one replaces a very
irregular boundary by a suitable polygonal or piecewise smooth
approximation then the eigenvalues are very little affected. This
continuous dependence of the spectrum on the boundary holds in
great generality for Dirichlet boundary conditions, but is much
less obvious for Neumann boundary conditions.

In the last part of the paper we investigate the effect of
perturbing the boundary. We first prove a quasi-monotonicity
property of the eigenvalues when the region decreases, under
suitable regularity hypotheses on the larger region. We then apply
a scaling trick to prove that the eigenvalues vary continuously
with the region provided the boundaries of the
regions concerned satisfy a uniform H\"older condition. Moreover
the change in the eigenvalues of the Neumann Laplacian is explicitly
estimated.

Many of the results of this paper apply not only to the Neumann
Laplacian but to general strictly elliptic second order operators
or \Schrodinger operators whose quadratic form domains are
contained in $W^{1,2}(\Ome )$. The proofs need almost no
alterations.

%%%%%%%%%%%%%%%%%%%%%%%%%%%%%%%%%%%%%%%%%%%%%%%%%%%%%%%%%%%%%%%%%%%%%%%%%
\section{Relationship between the Sobolev and Hardy-type Inequalities}

Let $\Omega\subset\R^N$ be a region with a finite measure
$|\Omega|$, for $x\in\Omega,$ $d(x)$ be the distance of the point
$x$ from the boundary $\partial\Omega$ of $\Omega$ and, for
$\varepsilon>0,$
$$
\partial_{\varepsilon}\Omega=\{x\in\Omega : d(x)\le\varepsilon\}.
$$
The {\it Minkowski dimension of} $~\partial\Omega$ {\it relative
to} $~\Omega$ (briefly, the {\it Minkowski dimension of}
$~\partial\Omega$) is the following quantity
$$
M(\partial\Omega)=\inf
{\{\lambda>0:M^{\lambda}(\partial\Omega)<\infty\}},
$$
where
$$
M^{\lambda}(\partial\Omega)= \limsup\limits_{\varepsilon\to0+}
\frac{|\partial_{\varepsilon}\Omega|}{\varepsilon^{N-\lambda}}.
$$
Obviously $M(\partial\Omega)\le N.$ However there exist $\Omega$
such that $M^{\lambda}(\partial\Omega)=\infty$ for all
$\lambda\in(0,N)$, \cite {EH}. It can be proved that
$M(\partial\Omega)\ge N-1$, \cite {Lap}. If $\Omega$ satisfies the
cone condition, then $M(\partial\Omega)=N-1$, \cite {Lap}.

Recall that a {\it Whitney covering} $\mathcal{W}$ of an open set
$\Omega$ is a family of closed cubes $Q$ each having edge length
$L_Q=2^{-k}, k=1,2,...,$ such that

(i) $\Omega=\bigcup\limits_{Q\in \mathcal{W}}Q;$

(ii) the interiors of distinct cubes are disjoint;

(iii) ${\rm
diam}\,(Q)\le{\rm dist}\,(Q,\partial\Omega)\le 4\,
{\rm diam}\,(Q);$

(iv) $\frac14{\rm diam}\,(Q_2)\le{\rm diam}\, (Q_1)\le 4\,{\rm diam}\,(Q_2)$
if $Q_1\bigcap Q_2\ne\O ;$

(v) at most $12^N$ other cubes in $\mathcal{W}$ can touch a fixed
$Q\in\mathcal{W}$, and for a fixed $t\in(1,5/4)$ each $x\in\Omega$
lies in at most $12^N$ of the dilated cubes $tQ$, $Q\in
\mathcal{W}.$

It is known  (see, for example, \cite [Chapter VI]{Ste}) that such
a covering exists for any $\Omega$. Note that condition (iii)
implies that ${\rm diam}\,(Q)\le d(x)\le5\,{\rm diam}\,(Q)$ for any
$x\in Q$.

Let, for a positive integer $k$, $n(k)$ denote the number of cubes
in $\mathcal{W}_k=\{Q\in \mathcal{W}:L_Q=2^{-k}\}$.
If
$\Omega$ has finite measure, then
$n(k)\le c_12^{Nk}$, where $c_1>0$ is independent of $k$.
Moreover, $M^{\lambda}(\partial\Omega)<\infty$ if, and only if, $n(k)\le c_2
2^{\lambda k}$, where $c_2>0$ is independent of $k$, \cite {MV}.

Let $0<\gamma\le1, M,\delta>0,s\ge1$ be an integer,
and let $\{V_j\}_{j=1}^s$ be a
family of bounded open cuboids and
$\{\lambda_j\}_{j=1}^s $ be a family of rotations. We say that,
for a bounded region $\Omega\subset \R^n$, its boundary $\partial\Omega\in
{\rm Lip}\,(\gamma,M,\delta,s,\{V_j\}_{j=1}^s,\{\lambda_j\}_{j=1}^s)$
if

(i) $ \Omega\subset \bigcup\limits_{j=1}^s(V_j)_{\delta}$, where
$(V_j)_{\delta}=\{x\in V_j: {\rm dist}\,(x,\partial
V_j)>\delta\},$ and $(V_j)_\delta\cap\Ome\ne\O;$

(ii) for $j=1,...,s$
$$
\lambda_j(V_j)=\{\,x\in \R^N:~a_{ij}<x_i<b_{ij}, \,i=1,....,N\},
$$
$$
\lambda_j(\Omega\cap V_j)=\{x\in\R^N:~a_{Nj}<x_N<\varphi_{j}(\bar x),~\bar
x\in W_j\},$$

where $\bar x=(x_1,...,x_{N-1})$, $W_j=\{\bar
x\in\R^{N-1}:~a_{ij}<x_i<b_{ij},\,i=1,...,N-1\}$
and
$$
|\varphi_j(\bar x)-\varphi_j(\bar y)|\le M\,|\bar x-\bar y|^\gamma, ~\bar
x, \,\bar y\in \overline W_j;
$$
(iii) if $V_j\cap\partial\Omega\ne\O$, then
$$
a_{Nj}+\delta\le\varphi_j(\bar x)\le b_{Nj}-\delta,~~\bar x\in
W_j.
$$
However, if $V_j\subset \Omega$, then $\varphi_j(\bar x)\equiv
b_{Nj}$.

We also say that, for a bounded region $\Omega$ and
$0<\gamma\le1,$ $\partial\Omega\in {\rm
Lip}\,\gamma$ if there exist $M,\delta>0$, an integer $s\ge1$,
a family of
bounded open cuboids $\{V_j\}_{j=1}^s$ and a family of
rotations $\{\lambda_j\}_{j=1}^s $ such that $\partial\Omega\in
{\rm Lip}\,(\gamma,M,\delta,s,
\{V_j\}_{j=1}^s,\{\lambda_j\}_{j=1}^s).$

If $\Omega$ is a bounded region and $\partial\Omega\in {\rm Lip}\,\gamma$,
then $M(\partial\Omega)\le N-\gamma$.
Moreover, if $\partial\Omega\in {\rm Lip}\,(\gamma,M,\delta,s,
\{V_j\}_{j=1}^s,\{\lambda_j\}_{j=1}^s)$,
there exist $\eps_0,a_0>0,$ depending only on $N,\gamma,M,\delta,s,
\{V_j\}_{j=1}^s,\{\lambda_j\}_{j=1}^s$, such that for all $0<\eps\le\eps_0$
\begin{equation}
|\partial_\eps\Ome|\le a_0\eps^\gamma.
\end{equation}
\begin{theorem}
Let $\Omega\subset\R^N$ be a region with a finite inradius, i.e.,
$\sup\limits_{x\in\Omega}d(x)<\infty$, and let $1\le p<\infty.$

$1$. If for some $\alpha>0, c_1>0$
\begin{equation}
\|d^{-\alpha}f\|_{L^p(\Omega)}\le c_1\,\|f\|_{W^{1,p}(\Omega)}
\end{equation}
for all $f\in W^{1,p}(\Omega)$, then
there exists $c_2>0$ such
that
\begin{equation}
\|f\|_{L^q(\Omega)}\le c_2\,\|f\|_{W^{1,p}(\Omega)}
\end{equation}
for all $f\in W^{1,p}(\Omega)$, where $q=\frac{Mp}{M-p},
M=N(1+\alpha)$ if $N>p$ and $q$ is any number such that
$p<q<p(1+\frac{\alpha p}N)$ if $N\le
p.$

$2$. If for some $\sigma>0$
\begin{equation}
\int\limits_\Omega d(x)^{-\sigma}{\rm d}^Nx<\infty
\end{equation}
and for some $q>p$
and $c_2>0$ inequality
$(4)$ holds, then there exists
$c_1>0$ such that inequality $(3)$ holds with $\alpha=\sigma (
\frac1p-\frac1q)$.
\end{theorem}

\Proof 1. First we note that there exists $c_3>0$ such that
$$
\|d^{N({1\over p}-{1\over r})}f\|_{L^r(\Omega)}
\le c_3\,\|f\|_{W^{1,p}(\Omega)}
$$
for all $f\in W^{1,p}(\Omega)$, where $r=\frac{Np}{N-p}$ if $N>p$ and
$p<r<\infty$ if $N\le p$. This inequality follows by scaling the standard
Sobolev inequality for cubes and by using the
Whitney decomposition of $\Omega$.
(It is contained in a more general statement of such type proved by R.C. Brown
\cite[Theorem 3.1]{Br}.)

Let $0<\lambda<1$ and $q\in(p,r)$ be such that
$$
{1\over q}={1-\lambda\over p}+{\lambda\over r}.
$$
Bearing in mind that
$$
|f|=
\bigg(d^{-N({1\over p}-{1\over r})
\frac{\lambda}{1-\lambda}}|f|\bigg)^{1-\lambda}
\bigg(d^{N({1\over p}-{1\over r})}
|f|\bigg)^{\lambda}
$$
we choose $\lambda$ so that
$N({1\over p}-{1\over r})
\frac{\lambda}{1-\lambda}=\alpha$.
Then
$$
\lambda=\frac{\alpha}{N}\bigg(\frac{\alpha}{N}+{1\over p}-{1\over r}\bigg)^{-1},
~q=\bigg(\frac{\alpha}{N}+{1\over p}-{1\over r}\bigg)
        \bigg[{1\over p}\bigg(\frac{\alpha}{N}+{1\over p}-{1\over r}\bigg)
-\frac{\alpha}{N}\bigg({1\over p}-{1\over r}\bigg)\bigg]^{-1}.
$$
By applying H\"older's inequality with the exponents $\frac{p}{1-\lambda}$
and ${r\over\lambda}$ we have
$$
\|f\|_{L^q(\Omega)}=\|(d^{-\alpha}|f|)^{1-\lambda}
(d^{N({1\over p}-{1\over r})}|f|)^{\lambda}\|_{L^q(\Omega)}
$$$$
\le\|d^{-\alpha}f\|_{L^p(\Omega)}^{1-\lambda}
\|d^{N({1\over p}-{1\over r})}f\|_{L^r(\Omega)}^{\lambda}
\le c_1^{1-\lambda}c_3^{\lambda}\|f\|_{W^{1,p}(\Omega)}.
$$
If $N>p$, then $r=\frac{Np}{N-p}$ and hence $q=\frac{Mp}{M-p}$. If
$N\le p$, then by passing to the limit as $r\to\infty$ we see that $q$
can be any real number satisfying $p<q<p(1+\frac{\alpha p}N)$.

2. The second statement follows immediately by H\"older's inequality with
the exponents $\frac{qp}{q-p}$ and $q$:
$$
\|d^{-\alpha}f\|_{L^p(\Omega)}\le\|d^{-\alpha}\|_{L^{\frac{qp}{q-p}}(\Omega)}
\|f\|_{L^q(\Omega)}
\le\bigg(\int\limits_\Omega d(x)^{-\sigma}{\rm d}^Nx\bigg)^{{1\over p}-{1\over
q}}c_2\,\|f\|_{W^{1,p}(\Omega)}.
$$
\begin{corollary}
Let $\Omega\subset\R^N$ be a region of finite measure and $1\le p<\infty$.
Then the following conditions are equivalent.

(a) For some $\alpha,c_1>0$ inequality $(3)$ holds for all
$f\in W^{1,p}(\Omega)$.

(b) For some $\sigma>0$ condition $(5)$ is satisfied
and for some $q>p$ and $c_2>0$ inequality $(4)$ holds for all
$f\in W^{1,p}(\Omega)$.

(c) $M(\partial\Omega)<N$
and for some $q>p$ and $c_2>0$ inequality $(4)$ holds for all
$f\in W^{1,p}(\Omega)$.
\end{corollary}

\Proof Inequality (3) implies, by putting
$f\equiv1$, that
$$
\int\limits_\Omega d(x)^{-\alpha p}{\rm d}^Nx\le c_1^p|\Omega|~.
$$
Now it suffices to recall that, for regions $\Omega$ of finite measure,
the inequality
$M(\partial\Omega)<N$ is equivalent to the existence of
$\mu\in(0,1)$ such that $\int_{\Omega}d(x)^{-\mu}\rmd^N x<\infty$,
\cite{Br}.

%%%%%%%%%%%%%%%%%%%%%%%%%%%%%%%%%%%%%%%%%%%%%%%%%%%%%%%%%%%%%%%%%%%%%%
\section{Equivalence of the Sobolev-type inequalities to
some spectral properties of Neumann Laplacian}

In this section we assume that $\Ome$ is a region
in $\R^N$ and suppose that $H=-\lap_N$
acts in $L^2(\Ome)$ subject to Neumann boundary conditions.

\begin{proposition}
Assume that $N\ge2$ and that $\Ome\subset\R^N$ is any region.

$1.$ Let $2<q\le {2N\over N-2}$ if $N\ge3$, $2<q<\infty$ if $N=2$ and
$M={2q\over q-2}$. Then the following conditions are equivalent.

(d) There exists $c_4>0$ such
that
\begin{equation}
\|f\|_{L^q(\Omega)}\le c_4\,\|f\|_{W^{1,2}(\Omega)} \label{smooth}
\end{equation}
for all $f\in W^{1,2}(\Omega)$.

(e) There exists $c_5>0$ such that
\begin{equation}
\norm \rme^{-Ht}f\norm_{L^\infty(\Ome)}
\leq c_5t^{-M/4}\norm f\norm_{L^2(\Ome)}
\end{equation}
for all $f\in L^2(\Ome)$ and all $0<t\leq 1$.

(f) The semigroup $\rme^{-Ht}$ has a continuous integral kernel
$K(t,x,y),t>0,x,y\in\Ome$ and there exists $c_6>0$ such that
\begin{equation}
0 <K(t,x,y)\leq c_6t^{-M/2} \label{kernel}
\end{equation}
for all $x,y\in \Ome$ and $0<t\leq 1$.

$2.$ Let $0<\gamma\le1$, $\Ome$ be bounded and $\partial\Ome\in{\rm Lip}\,
\gamma.$ If $\gamma=1,$ then (d) is satisfied with $q={2N\over N-2}$ for
$N\ge3$ $($hence in $(7)$ and $(8)$ $M=N$$)$ and with
any $2<q<\infty$ for $N=2$
$($hence in $(7)$ and $(8)$ any $M>2$$)$. If $0<\gamma<1,$
then (d) is satisfied
with $q={2(\gamma+n-1)\over N-1-\gamma}$ $($hence in $(7)$ and $(8)$
$M={\gamma+N-1\over\gamma}$$)$.
\end{proposition}

The first statement is proved, for example, in \cite[Corollary 2.4.3, Lemma
2.1.2]{HKST}. (One needs to take into account that
${\rm Quad}\,(H)=W^{1,2}(\Ome).$) The second statement is proved in
\cite{Glo}, \cite{Maz1}, \cite{Maz2}.

\begin{remark}
Each of the constants $c_4,c_5,c_6$ can be estimated from any of the others,
given $q$, or equivalently $M$.
\end{remark}
\begin{remark}
If $N=1$, then ({\it d\,}) is satisfied with $q=\infty$ and ({\it e\,}) and
({\it f\,}) are satisfied
with $M=1$.

If $N>1$, there exists a
region, say a region  with exponentially degenerate boundary \cite{Maz2},
such that ({\it d\,}) is not valid for any $q>2$. However, if such $q>2$
exists, it must satisfy the assumptions of  Proposition 3.
The appropriate range for $M={2q\over q-2}$ is $N\le M<\infty$
for $N\ge3$ and $2<M<\infty$ for $N=2$.

On the other hand ({\it e\,}) and ({\it f\,}), which are always equivalent
\cite[Lemma 2.1.2]{HKST},
could be also invalid for some region $\Omega$ for all $M>0.$ However,
if there exist $M>0$ for which ({\it e\,}) and ({\it f\,}) are valid,
then $N\le M<\infty$ for any $N\ge1$.

Hence, if $N\le M<\infty$ for $N\ge3$ and  $2<M<\infty$ for $N=2$, then
({\it e\,}) and ({\it f\,}) are equivalent to ({\it d\,})
where $q={2M\over M-2}.$ However, if
$M=2$ for $N=2$, then ({\it e\,}) and ({\it f\,}) are not equivalent
to ({\it d\,})
for any $q>2$. (In this case ({\it e\,}) and ({\it f\,}) are equivalent to a
certain
logarithmic Sobolev inequality or to a certain Nash inequality
\cite[Example 2.3.1, Corollary 2.4.7]{HKST}.)
\end{remark}

\begin{example}
Let $N\ge3$ and $0<\gamma\le1$ or $N=2$ and $0<\gamma<1$. The following well
known example shows that in this case, the exponents $q$ and $M$ in the
second statement of Proposition 3 are the best possible, i.e., $q$ cannot be
replaced by a larger one and $M$ cannot be replaced by a smaller one.
Let $\Ome=\{(x,y):y\in\R^{N-1},|y|<1, |y|^\gamma<x<1\}.$ Then $\partial\Ome\in
{\rm Lip}\,\gamma.$ A direct computation shows that
$x^{-\delta}\in W^{1,2}(\Ome)$ if, and only if, $\delta<
-1+{1\over2}(1+{N-1\over\gamma})$ and
$x^{-\delta}\in L^q(\Ome)$ if, and only if,
$\delta<{1\over q}(1+{N-1\over\gamma})$.
If ({\it d}\,) holds, then
$${1\over q}\bigg(1+{N-1\over\gamma}\bigg)
\ge-1+{1\over2}\bigg(1+{N-1\over\gamma}\bigg)\Longleftrightarrow
q\le\frac{2(\gamma+N-1)}{N-1-\gamma}.$$
Since in the case under consideration
({\it e}\,) and ({\it f}\,) are equivalent to ({\it d}\,) it follows also that
$M\ge\frac{\gamma+N-1}\gamma.$
\end{example}
\begin{theorem}
Assume that $\Ome\subset\R^N$ is a region of finite measure.

$1.$ The following conditions are equivalent.

(g) For some $q>2$ and $c_4>0$ the inequality
$$
\|f\|_{L^q(\Omega)}\le c_4\,\|f\|_{W^{1,2}(\Omega)}
$$
is satisfied for all $f\in W^{1,2}(\Omega)$.

(h) $H$ has discrete spectrum and if all its eigenvalues $\lam_n, n=0,1,2...$,
which are nonnegative and of finite multiplicity, are
written in increasing order and repeated according to multiplicity
and $f_n$ is the corresponding orthonormal basis of eigenvectors,
then there exist $\alpha_1,c_7,\alp_2,c_8>0$ and an integer $n_0\ge1$ such that
\begin{equation}
\lam_n\geq c_7n^{\alp_1},~~~\norm f_n
\norm_{L^\infty(\Ome)}
\leq c_8\lam_n^{\alp_2}
\end{equation}
for all $n\ge n_0$.

$2.$ If $N=1$ or $N\ge2$, $\Ome$ is bounded and $\partial\Ome\in
{\rm Lip}\,\gamma$ where $0<\gamma\le1,$ then conditions (g) and (h) are
satisfied.
\end{theorem}
\begin{remark}
One could also assume that conditions (9) were valid for all $n\ge1$.
However, the first few eigenvalues may be extremely small if $\Omega$ is
nearly  disconnected and it is not easy to provide explicit bounds on the
constants $c_7,c_8$ which apply for all $n\ge1$.
\end{remark}
\begin{remark}
The relationship between $q$ and $(\alpha_1,\alpha_2)$ is not
symmetrical and we do not expect that a symmetrical
relationship can be obtained.
\end{remark}

The proof of this theorem will be based on the following lemmas containing
additional information.

\begin{lemma}
Let $M,c_5>0$ and let $\Ome\subset\R^N$ be a region of finite
measure such that
inequality $(7)$ is satisfied for all $0<t\le1$. Then
\begin{equation}
\|f_n\|_{L^\infty(\Ome)}\le c_{9}\cases{1 &~ {\rm if} $0<\lam_n\le1$,\cr
\lam_n^{M\over4} &~ {\rm if}~~~~~~ $\lam_n>1$,\cr}
\end{equation}
where $c_{9}={\rm e}c_5.$
\end{lemma}
\Proof
If $0<\lam_n\le1$ put $t=1$ in (7)  to get
${\rm e}^{-\lam_n}\|f_n\|_{L^\infty(\Ome)}
\le c_5.$ So $\|f_n\|_{L^\infty(\Ome)}\le c_5{\rm e}^{\lam_n}\le c_5{\rm e}.$
If $\lam_n>1$ put $t=1/\lam_n$ in (7)  to get
${\rm e}^{-1}\|f_n\|_{L^\infty(\Ome)}\le
c_5\lam_n^{M\over4}.$

\begin{lemma}
Let $M,c_6, c_{10}>0$ and let $\Ome\subset\R^N$ be a region
such that inequality $(8)$ is satisfied for all $x,y\in\Ome,0<t\le1$ and
$|\Ome|\le c_{10}.$ Then
there exists an integer $n_0\ge1,$ depending only on $c_6,c_{10},$ such that
\begin{equation}
\lam_n\ge\bigg(\frac{n}{n_0}\bigg)^{2\over M},~~n\ge n_0.
\end{equation}
\end{lemma}
\Proof
By integrating (8) with $x=y$ over $\Ome$ we get
$$
n\rme^{-\lam_nt}\le\sum\limits_{k=0}^n\rme^{-\lam_kt}
\le\sum\limits_{k=0}^\infty\rme^{-\lam_kt}
=\int_\Ome K(t,x,x)\rmd^Nx\le c_6|\Ome|t^{-{M\over2}},
$$
hence
$$
\rme^{\lam_nt}\ge\frac{nt^{M\over2}}{c_6|\Ome|}
\ge\frac{nt^{M\over2}}{c_6c_{10}}, ~~~0<t\le1.
$$
By putting here $t=1$ it follows that
$$
\lam_n\ge1, ~~~ n\ge n_0=[\rme c_6c_{10}]+1.
$$
Finally, for $n\ge n_0$ we put $t=\lam_n^{-1}$ to get
$$
\lam_n\ge\bigg(\frac n{\rme c_6c_{10}}\bigg)^{2\over M}\ge
\bigg(\frac n{n_0}\bigg)^{2\over M}.
$$

\begin{lemma}
Let $M,c_{9}>0$ and let $n_0\ge1$ be an integer. Moreover, let
$\Ome\subset\R^N$ be a region of finite measure such that
inequalities $(10)$ and $(11)$ are satisfied. Then there there
exist $c_5,c_6>0,$ depending only on $M, c_9$ and $n_0,$ such that
the inequalities $(7)$, $(8)$ are satisfied with $2M$ replacing
$M$.
\end{lemma}
\Proof
If $0<t\le1$ and $x,y\in\Ome,$ then
$$
0<K(t,x,y)=\sum\limits_{n=0}^\infty\rme^{-\lam_nt}f_n(x)f_n(y)
\le\sum\limits_{n=0}^\infty\rme^{-\lam_nt}\|f_n\|_{L^\infty(\Ome)}^2
$$$$
\le c_9^2\bigg(\sum\limits_{n=0}^{n_0-1}1+
\sum\limits_{n=n_0}^\infty\rme^{-\lam_nt}\lam_n^{M\over2}\bigg).
$$
Since
$$
\rme^{-\lam_nt/2}\lam_n^{M\over2}
=\bigg({t\over2}\bigg)^{-{M\over2}}\rme^{-{{\lam_nt}/2}}
\bigg({{\lam_nt}\over2}\bigg)^{M\over2}
\le\bigg({t\over2}\bigg)^{-{M\over2}}
\bigg({M\over2\rme}\bigg)^{M\over2}
=c_{11}t^{-{M\over2}}
$$
it follows that
$$
\sum\limits_{n=n_0}^\infty\rme^{-\lam_nt}\lam_n^{M\over2}
\le c_{11}t^{-{M\over2}}
\sum\limits_{n=n_0}^\infty\rme^{-{\lam_nt\over2}}
$$$$
\le c_{11}n_0t^{-{M\over2}}{1\over n_0}
\sum\limits_{n=n_0}^\infty
\rme^{-{1\over2}({n\over n_0})^{2\over M}t}
\le c_{11}n_0t^{-{M\over2}}
\int_0^\infty
\rme^{-{1\over2}x^{2\over M}t}\rmd x
$$$$
=c_{11}n_0t^{-{M}}
\int_0^\infty
\rme^{-{1\over2}s^{2\over M}}\rmd s
=c_{12}t^{-M}.
$$
Hence
\begin{equation}
\sum\limits_{n=0}^\infty\rme^{-\lam_nt}\|f_n\|_{L^\infty(\Ome)}^2
\le c_6t^{-M}, \label{S}
\end{equation}
where $c_6=c_9^2(n_0+c_{12}),$ and (8) follows with the exponent
$-{M\over2}$ replaced by $-M$.

Furthermore,
$$
\|\rme^{-Ht}f\|_{L^\infty(\Ome)}
=\|\sum\limits_{n=0}^\infty\rme^{-\lam_nt}f_n(x)(f,f_n)\|_{L^\infty(\Ome)}
$$$$
\le\sum\limits_{n=0}^\infty\rme^{-\lam_nt}\|f_n\|_{L^\infty(\Ome)}|(f,f_n)|
\le\bigg(\sum\limits_{n=0}^\infty
\rme^{-2\lam_nt}\|f_n\|_{L^\infty(\Ome)}^2\bigg)^{1\over2}
\|f\|_{L^2(\Ome)}.
$$
Hence by (\ref{S})
$$
\|\rme^{-Ht}f\|_{L^\infty(\Ome)}
\le (2^{-M}c_6)^{1\over2}t^{-{M\over2}}\|f\|_{L^2(\Ome)}
$$
and (7) follows with $c_5=(2^{-M}c_6)^{1\over2}$ and the exponent
$-{M\over4}$ replaced by $-{M\over2}$.
\vskip0.3cm
\noindent
\underline{Proof of Theorem 7} By Proposition 3 and Lemmas 10, 11
$(g)$ implies $(h)$ with $\alpha_1={2\over M}$ and $\alpha_2=
{M\over4}$ where $M={2q\over q-2}.$ Conversely by the proof of
Lemma 12 it follows that $(h)$ implies inequality (7) with
$M=2(2\alpha_2+{1\over\alpha_1})$ which in its turn by Proposition 3
implies $(g)$ with $q={2M\over M-2}$.

%%%%%%%%%%%%%%%%%%%%%%%%%%%%%%%%%%%%%%%%%%%%%%%%%%%%%%%%%%%%%%%%%%%%%%%%%%%
\section{Perturbations of the Domain}

In this section we compare the spectrum of $H_i=-\lap_N$ acting in
$L^2(\Ome_i)$ when $\Ome_1$ and $\Ome_2$ are very close to each
other in a suitable sense. We will also need to assume regularity,
since it is known that even if $\Ome_1$ has smooth boundary and
$\Ome_2$ only differs from it in an arbitrarily small neighbourhood
of a single point of $\partial\Ome_1$, the spectrum of $H_2$ need
not be discrete.

We start with the more general argument. Following
\cite[Chapter 4]{STDO} we define the variational quantities $\mu_{n,i}$
for all non-negative integers $n$ by
\[
\mu_{n,i}=\inf \{ \mu(L):\dim (L)=n+1\}
\]
where $\mu(L)$ is defined for every finite-dimensional subspace $L$
of $L^2(\Ome_i )$ by
\[
\mu(L)=\sup\{ Q_i(f)/\norm f \norm_{L^2(\Ome)}^2 :0\not= f\in L\}
\]
and $Q_i$ are defined as in (\ref{B}). Note that $\mu_{0,i}=0$
since $0$ is an eigenvalue of $H_i$ for $i=1,2$. It is known that
$\mu_{n,i}$ are equal to the eigenvalues $\lam_{n,i}$ of $H_i$
written in increasing order and repeated according to multiplicity
in case $H_i$ has compact resolvent.

\begin{theorem}\label{M}
Let $\Ome_1\subset\R^N$ be a region of finite measure.

$1.$ If for some $q>2, c_{13}>0$
\begin{equation}
\|f\|_{L^q(\Omega_1)}\le c_{13}\,\|f\|_{W^{1,2}(\Omega_1)} \label{C0}
\end{equation}
for all $f\in W^{1,2}(\Omega_1)$, then
for every integer $n\geq 1$
there exist $b_{n,1}=b_{n,1}(\Ome_1),\eps_{n,1}=\eps_{n,1}(\Ome_1)>0$
such that for all regions
$\Ome_2\subset\Ome_1$, satisfying $|\Ome_1\setminus\Ome_2|\le \eps_{n,1}$,
the inequality
\begin{equation}
\mu_{n,2} \leq (1+b_{n,1}|\Ome_1\setminus\Ome_2|)\lam_{n,1}   \label{C1}
\end{equation}
holds.

$2.$ If, in addition, $M(\partial\Ome_1)<N,$
then for every $\sigma\in(0,N-M(\partial
\Ome_1))$ and for every integer $n\geq 1$
there exist $b_{n,2}=b_{n,2}(\Ome_1),\eps_{n,2}=\eps_{n,2}(\Ome_1)>0$
such that
for all $0<\eps\le\eps_{n,2}$ and for all regions $\Ome_2$, satisfying
$\Ome_1\setminus\partial_\eps\Ome_1\subset
\Ome_2\subset \Ome_1$, the inequality
\begin{equation}
\mu_{n,2} \leq (1+b_{n,2}\eps^\sig)\lam_{n,1}   \label{C2}
\end{equation}
holds.
\end{theorem}

\Proof $1.$ Let $L=\lin\{ \phi_0,...,\phi_n\}$
where $\phi_i$ are the eigenfunctions of $H_1$ associated with
$\lam_{i,1}$ and let $M=PL$ where $P$ is the restriction map from
$L^2(\Ome_1)$ to $L^2(\Ome_2)$.
By Proposition 3 it follows that there exist $M>2$
and $c_5>0$ such that inequality (7) is satisfied with $\Ome_1$ and $H_1$
replacing $\Ome$ and $H$ for all $0<t\le1$ and $f\in L^2(\Ome_1)$. Hence if
$f=\sum\limits_{k=0}^n\alpha_k\phi_k\in L$, $\norm f\norm_{L^2(\Ome)}=1$,
then by (7) where $t=1$ applied to $\rme^{H_1}f$ we have
$$
\norm f\norm_{L^\infty(\Ome)}\leq c_5\norm \rme^{H_1}  f\norm_{L^2(\Ome)}
=c_5\,\norm\sum\limits_{k=0}^n\alpha_k\rme^{\lam_{k,1}}\phi_k\norm_{L^2(\Ome)}
$$$$
\leq
c_5\,(\sum\limits_{k=0}^n|\alpha_k|^2\rme^{2\lam_{k,1}})^{1\over2}
\le c_5\,\rme^{\lam_{n,1}}.
$$
Furthermore
$$
\norm Pf \norm_{L^2(\Ome_1)}^2=\norm f\norm_{L^2(\Ome_2)}^2=\norm f\norm_{L^2(\Ome_1)}^2
-\norm f\norm_{L^2(\Ome_1\setminus\Ome_2)}^2
$$$$
\ge 1-|\Ome_1\setminus\Ome_2|
\,\norm f\norm_{L^\infty(\Ome_1\setminus\Ome_2)}^2
\ge 1-c_5^2\,\rme^{2\lam_{n,1}}|\Ome_1\setminus\Ome_2|.
$$

Assume that $|\Ome_1\setminus\Ome_2|\le(2c_5^2{\rm e}^{2\lam_{n,1}})^{-1}$.
Then
$$
\norm Pf \norm_{L^2(\Ome_1)}^{-2}
\le1+2c_5^2\,\rme^{2\lam_{n,1}}|\Ome_1\setminus\Ome_2|.
$$

If $g=Pf$ where $f\in L$ and $\norm f \norm_{L^2(\Ome_1)}=1$, then
$$
\frac{Q_2(g)}{\norm g\norm_{L^2(\Ome_2)}^2}
=\frac{Q_2(f)}{\norm
Pf\norm_{L^2(\Ome_1)}^2}
\leq (1+2c_5^2\,\rme^{2\lam_{n,1}}|\Ome_1\setminus\Ome_2|)
Q_1(f)
$$

Since $\dim(M)\le\dim(L)=n+1$, the first statement of the theorem
with $b_{n,1}=2c_5^2{\rm e}^{2\lam_{n,1}}$
and $\eps_{n,1}=b_{n,1}^{-1}$ follows
by using the variational definitions of $\mu_{n,2}$ and $\lam_{n,1}$.

$2.$ If $M(\partial\Ome_1)<N$, then for every
$\sigma\in(0,N-M(\partial\Ome))$ there exist $\eps_1=\eps_1(\sigma),
a_1=a_1(\sigma)>0$
such that for all $0<\eps\le\eps_1$
\begin{equation}
|\partial_\eps\Ome_1|\le a_1\eps^\sigma,
\label{C3}
\end{equation}
Hence the second statement of the theorem with $b_{n,2}=b_{n,1}a_1$ and
$0<\eps\le \eps_{n,2}\equiv {\rm min}\,\{b_{n,2}^{-{1\over\sig}},\eps_1\}$
follows from the first one.

\begin{remark}
The size of $b_{n,1}(\Ome)$ and $\eps_{n,1}(\Ome)$ depends upon
$\lam_{n,1}$. An upper bound to $\lam_{n,1}$ can be given in terms
of the inradius $r=\,{\rm max}\,\{d(x):x\in\Ome\}$ as follows. If
$B(a,r)\subset\Ome,$ then $\lam_{n,1}\le\gamma_{n,a,r}$ where
$\gamma_{n,a,r}$ is the {\it n}\,th eigenvalue of $-\Delta$ in
$L^2(B(a,r))$ subject to Dirichlet boundary conditions. By scaling
one also has $\gamma_{n,a,r} =\gamma_{n,0,1}r^{-2}$.
\end{remark}
\begin{remark}
If for some $\alpha,c_{14}>0$
$$
\|d^{-\alpha}f\|_{L^2(\Omega_1)}\le c_{14}\,\|f\|_{W^{1,2}(\Omega_1)}
$$
for all $f\in W^{1,2}(\Omega_1)$, then both statements of Theorem 10
are valid by Corollary 2.
\end{remark}

The conditions of Theorem \ref{M} are not sufficient to establish
that $H_2$ has a compact resolvent, since $\partial\Ome_2$ may have
arbitrarily bad local singularities subject to the above
conditions. In order to obtain an inequality in the reverse
direction we make further assumptions.

\begin{corollary}
Assume that $\Ome_1$ satisfies the conditions of the first part of
Theorem $10$ and for some $\sigma>0$ inequality $(\ref{C3})$ holds. Moreover,
let regions $\Ome_3(\eps),\eps>0,$ be such that

~~$($i$\,)$~~ $\Ome_3(\eps)\subset\Ome_1\setminus\partial_\eps\Ome_1,$

~$($ii$\,)$~~there exist $\eps_2,a_2>0$ such that for all $0<\eps\le\eps_2$
$$
|\Ome_1\setminus\Ome_3(\eps)|\le a_2\eps^\sigma,
$$

$($iii$\,)$~~for every integer $n\ge1$ there exist $b_{n,3}=b_{n,3}(\Ome_1),
\eps_{n,3}=\eps_{n,3}(\Ome_1)>0$
such that for all $0<\eps\le\eps_{n,3}$
$$
\mu_{n,3}\geq \lam_{n,1}(1- b_{n,3}\eps^\sig).\label{E}
$$
Then for every integer $n\ge1$ there exist $b_{n,4}=b_{n,4}(\Ome_1),
\eps_{n,4}=\eps_{n,4}(\Ome_1)>0$
such that
for all $0<\eps\le\eps_{n,4}$ and
for every region $\Ome_2$, for which
inequalities $(\ref{C0})$ and $(\ref{C3})$ holds with
$\Ome_2$ replacing $\Ome_1$
$($with the same $q,c_{13},\sig,a_1,\eps_1$$)$ and
$\Ome_1\setminus\partial_\eps\Ome_1\subset\Ome_2\subset\Ome_1$, the inequality
\begin{equation}
(1-b_{n,4}\eps^\sig)\lam_{n,1}\le\lam_{n,2} \leq (1+b_{n,4}\eps^\sig)\lam_{n,1}
\label{D}
\end{equation}
holds.
\end{corollary}
\Proof An application of Theorem \ref{M} to the pair $\Ome_2,\Ome_1$ yields
\[
\lam_{n,2}\leq
(1+b_{n,2}\eps^\sig)\lam_{n,1}.
\]
for $0<\eps\le\eps_{n,2}$. In particular,
\[
\lam_{n,2}\leq
(1+b_{n,2}\eps_{n,2}^\sig)\lam_{n,1}.
\]
An application of Theorem \ref{M} to the pair
$\Ome_3(\eps),\Ome_2$ yields
\[
\mu_{n,3}\leq
(1+b_{n,5}\eps^\sig)\lam_{n,2}
\]
for $0<\eps\le\eps_{n,5}$, where $b_{n,5}=2c_5^2{\rm e}^{2\lam_{n,2}}$ and
$\eps_{n,5}={\rm min}\,\{b_{n,5}^{-{1\over\sig}},\eps_1\}$. Note that
$b_{n,5}\le b_{n,6}$ and $\eps_{n,5}\ge \eps_{n,6}$ where
$b_{n,6}=2c_5^2{\rm exp}\,(2(1+b_{n,2}\eps_{n,2}^\sig)\lam_{n,1})$ and
$\eps_{n,2}={\rm min}\,\{b_{n,6}^{-{1\over\sig}},\eps_1\}$
depend only on $n,\sig$ and $\Ome_1$. Together with assumption ({\it iii})
this yields (\ref{D}).

The above theorem may be applied to regions with
Lip$\,\gamma$ boundaries.
We start with the simplest
example. Let $0<\gamma\le1,M,k>0$ and let
\[
\Ome_i=\{x\in\R^N:0<x_N<\phi_i(\bar x),\bar x\in G\}, ~i=1,2,
\]
where $G$ is a bounded region in $\R^{N-1}$ with a smooth boundary.
We assume that

\[
|\phi_i(\bar x)-\phi_i(\bar y)|\leq M|\bar x-\bar y|^\gam,~~\bar x,\bar y
\in \bar G,~~i=1,2,
\]
and
\[
k^{-1}\leq \phi_i(\bar x)\leq k,~~\bar x\in \bar G,~~i=1,2.
\]
We do not assume any relationship between the directions of
normals of $\Ome_1$ and $\Ome_2$, or even that these normal
directions exist.

\begin{lemma}\label{P}
Under the conditions of the last paragraph for every integer $n\ge1$
there exist $b_{n,7}=b_{n,7}(\Ome_1),
\eps_{n,7}=\eps_{n,7}(\Ome_1)>0$ such that for all $0<\eps\le\eps_{n,7}$
and all $\phi_2$, satisfying
$$
(1-\eps)\phi_1(\bar x)\le\phi_2(\bar x)\le\phi_1(\bar x),~~\bar x\in\bar G,
$$
the inequality
\[
(1- b_{n,7}\eps)\lam_{n,1}\leq \lam_{n,2}\leq
(1+b_{n,7}\eps)\lam_{n,1}
\]
holds.
\end{lemma}

\Proof
Since $\partial\Ome_1\in {\rm Lip}\,\gamma$ it follows, as noted in
Proposition 3, that inequality (\ref{C0}) is valid for some $q,c_{13}>0$.
Moreover it also holds with $\Ome_2$ replacing $\Ome_1$ (with the same
$q,c_{13}>0$). Since the operator $H_2$ has compact resolvent and
$|\Ome_1\setminus\Ome_2|\le k\eps|G|$, inequality (\ref{C1}) yields
$$
\lam_{n,2}\le(1+b_{n,1}k|G|\eps)\lam_{n,1}.
$$
Similarly for
$$
\Ome_3(\eps)=\{x\in\R^N:0<x_N<(1-\eps)\phi_1(\bar x),\bar x\in G\}
$$
where $0<\eps<1/2$ we have
$$
\lam_{n,3}\le(1+b_{n,1}k|G|\eps)\lam_{n,2}.
$$

To derive the estimate below for $\lam_{n,3}$ we transfer the
quadratic form $Q_3$ to $L^2(\Ome_1)$ by means of the unitary map
$U_\eps:L^2(\Ome_1)\to L^2(\Ome_3(\eps))$ defined by
\[
(U_\eps f)(\bar x,x_N)=(1-\eps)^{-1/2}f(\bar x,(1-\eps)^{-1}x_N).
\]
The inequality
\[
Q_1(f)\leq Q_3(U_\eps f)
\]
valid for all $f\in W^{1,2}(\Ome_1)$, yields the inequality
$\lam_{n,1}\leq\lam_{n,3}$ by the variational method.
\vskip0.3cm
We now turn to the application of Theorem \ref{M} to a general
region of H\"older type. The proof of our main result, Theorem 21,
depends upon the construction of mappings $T_\eps$ of $\Ome$ into
itself satisfying the properties (20), (22) and (25) below.

Other
definitions of regions of H\"older type are possible but Theorem
20 is still valid for such definitions provided similar mappings
can be constructed.
The underlying idea of that theorem can also be applied to
uniformly elliptic operators of the form
\[
Hf=-\sum_{i,j=1}^N\frac{\partial}{\partial x_i}\left\{
a_{i,j}(x) \frac{\partial f}{\partial x_j}\right\}
\]
subject to Neumann boundary conditions provided the coefficients
are H\"older continuous in some neighbourhood of the boundary.

Let a bounded region $\Ome\subset\R^N$ be such that
$\partial\Omega\in
{\rm Lip}\,(\gamma,M,\delta,s, \{V_j\}_{j=1}^s,\{\lambda_j\}_{j=1}^s)$.
Then also
$\partial\Omega\in
{\rm Lip}\,(\gamma,M,{\delta\over4},s, \{(V_j)_{\delta\over2}\}_{j=1}^s,
\{\lambda_j\}_{j=1}^s)$.
Note that $(V_j)_{\delta\over2}=\lambda_j^{-1}((W_j)_{\delta\over2}\times
(\tilde a_{Nj}, \tilde b_{Nj}))$ where
$\tilde a_{Nj}=a_{Nj}+{\delta\over2}, \tilde b_{Nj}=b_{Nj}-{\delta\over2},$
and, in addition to conditions 1) and 2) of the appropriate definition, also
the following condition is satisfied
\begin{equation}\label{ZZ}
\partial\Ome\cap\lam_j^{-1}((W_j)_{\delta\over2}\times[\tilde a_{Nj}
-{\delta\over4},\tilde a_{Nj}])=\O,~~j=1,...s.
\end{equation}

Moreover, let functions $\psi_j\in C^{\infty}(\R^N)$ satisfy
$0\le\psi_j\le1,{\rm supp}\,\psi_j\subset (V_j)_{{3\over4}\delta},
|\nabla\psi_j|\le {b\over\delta},$ where $b>0$ is a constant,
$j=1,...,s,$ and $\sum\limits_{j=1}^s\psi_j(x)=1$ for $x\in\Ome$.

Let $e_N=(0,...,0,1)$ and $\xi_j=\lam_j^{-1}(e_N), j=1,...,s.$ For
$x\in\R^N$ and $\eps\in(0,{\delta\over4}]$ define
\begin{equation}
T_\eps(x)=x-\eps\sum\limits_{j=1}^s\xi_j\psi_j(x).\label{T}
\end{equation}
\begin{lemma}
Let a bounded region $\Ome\subset\R^N$ be such that
$\partial\Omega\in
{\rm Lip}\,(\gamma,M,\delta,s,$
$\{V_j\}_{j=1}^s,\{\lambda_j\}_{j=1}^s)$.
Then there exist $A_1,A_2,A_3,\eps_3>0,$ depending only on
$N,\gamma,M,$
$\delta,s,\{V_j\}_{j=1}^s,
\{\lambda_j\}_{j=1}^s,$
such that for all $0<\eps\le\eps_3$
\begin{equation}\label{A0}
\bigg|\frac{\partial T_{\eps i}}{\partial
x_j}(x)-\delta_{ij}\bigg|\le A_1\eps,~~x\in\R^N,
\end{equation}
in particular,
the Jacobian determinant ${\rm Jac}\,(T_\eps,x)$
satisfies the inequality
\begin{equation}\label{A3}
{1\over2}\le 1-A_2\eps\le{\rm Jac}\,(T_\eps,x)\le1+A_2\eps,~~x\in\R^N.
\end{equation}
Moreover,
\begin{equation}\label{A1}
\Ome\setminus\partial_\eps\Ome\subset T_\eps(\Ome)\subset\Ome
\setminus\partial_{A_3\eps^{1\over\gamma}}\Ome.
\end{equation}
\end{lemma}
\Proof
1. Since the Jacobi matrix of the map $T_\eps$ has the form $I+\eps
B$ where $I$ is the identity matrix and $B$ is a matrix whose
elements $b_{ij}$ are independent of $\eps$ and bounded:
$|b_{ij}|\le {bs\over\delta}$, it follows that there exists
$A_1>0$, depending only on $N,\delta$ and $s$ such that inequality
(\ref{A0}) is satisfied for all $0<\eps\le1$. Hence, for all
sufficiently small $\eps>0$ and for all $x\in\R^N$ inequality
(\ref{A3}) is satisfied. Consequently, for all those $\eps$ the map
$T_\eps:\R^N
\to\R^N$ is one-to-one. Indeed, it is locally one-to-one since
${\rm Jac}\,(T_\eps,x)\ge{1\over2}$ and it is also globally
one-to-one since $|x-y|>2\eps$ implies $T_\eps(x)\ne T_\eps(y).$
Also $T_\eps(\Ome)$ is a region and
$T_\eps(\partial\Ome)=\partial T_\eps(\Ome).$

2. For $x\in\R^N$ let
$$J(x)=\{j\in\{1,...,s\}:
x\in (V_j)_{{3\over4}\delta}\}.$$
The inclusion ${\rm supp}\,\psi_j\subset
(V_j)_{{ 3\over4}\delta}$ implies that $\psi_j(x)=0$ for $j\not\in J(x)$ and
$$
T_\eps(x)=x-\eps\sum\limits_{j\in J(x)}\xi_j\psi_j(x).
$$
Let
$$
C(x)=\{\sum\limits_{j\in J(x)}\alpha_j\xi_j,~\alpha_j>0\}
$$
and
$$
C(x,\eps)=\{\eps\sum\limits_{j\in J(x)}\alpha_j\xi_j,~\alpha_j>0,
\sum\limits_{j\in J(x)}\alpha_j<1\}.
$$
We claim that
\begin{equation}\label{A4}
(x-C(x))
\cap(\bigcap\limits_{j\in J(x)}V_j)\subset\Ome,~~x\in\bar\Ome,
\end{equation}
and
\begin{equation}\label{A4'}
(x-C(x,\eps))
\subset\Ome\cap(\bigcap\limits_{j\in
J(x)}V_j),~~x\in\bar\Ome,0<\eps\le{\delta \over4}.
\end{equation}
Also
\begin{equation}\label{A4''}
(x-\overline{C(x)})
\cap(\bigcap\limits_{j\in J(x)}V_j)\subset\Ome,~~x\in\Ome,
\end{equation}
and
\begin{equation}\label{A4'''}
(x-\overline{C(x,\eps)})
\subset\Ome\cap(\bigcap\limits_{j\in
J(x)}V_j),~~x\in\Ome,0<\eps\le{\delta \over4},
\end{equation}
which implies, in particular, that $T_\eps(\Ome)\subset\Ome.$

        Indeed, let $x\in\bar\Ome$ and $j_1\in J(x)$. Since $x\in
V_{j_1}\bigcap\Ome$  and
$\lam_{j_1}(V_{j_1}\cap\Ome)$ is a subgraph
it
follows that $\{x-\alpha_{j_1}\xi_{j_1},
\alpha_{j_1}>0\}\cap V_{j_1}\subset\Ome$. Next, let
$j_2\in J(x),j_2\ne j_1,$ and for some $\alpha_{j_1}>0$,
$x-\alpha_{j_1}\xi_{j_1}\in V_{j_1}\cap V_{j_2}.$
For similar reasons it follows that
$
\{x-\alpha_{j_1}\xi_{j_1}-\alpha_{j_2}\xi_{j_2}, \alpha_{j_2}>0\}\cap
V_{j_2} \subset\Ome,
$
hence
$$
\{x-\alpha_{j_1}\xi_{j_1}-\alpha_{j_2}\xi_{j_2},
\alpha_{j_1},\alpha_{j_2}>0\} \cap V_{j_1}\cap
V_{j_2} \subset\Ome
$$
and so on. Since for $y\in C(x,\eps)$
$$
|y|=\eps\sum\limits_{j\in J(x)}\alpha_j\xi_j\le
\eps\sum\limits_{j\in J(x)}\alpha_j\le\eps\le{\delta\over4},
$$
condition (\ref{ZZ}) implies that for all $j\in J(x)$
$$
x-C(x,\eps)\subset ((V_j)_{3\delta\over4})^{\delta\over4}\subset V_j.
$$
Hence
$
x-C(x,\eps)\subset\bigcap\limits_{j\in J(x)}V_j
$
and, by (\ref{A4}), (\ref{A4'}) follows. If $x\in\Ome,$ then in the
argument above one may assume that $\alpha_j\ge0,j\in J(x)$, hence
(\ref{A4''}) and (\ref{A4'''}) follow.
        %We claim  that for $x\in\bar\Ome$ and $\eps\in(0,{\delta\over4}]$
%$$
%S(x)=\{x-\eps\sum\limits_{j\in J(x)}\alpha_j\xi_j,~0<\alpha_j\le1\}
%\subset\Ome,
%$$
%which means, in particular, that $T_\eps(\Ome)\subset\Ome.$
%(If $x\in\Ome$ then $\overline{S(x)}\subset\Ome.$)
%
%Indeed, let $j_1\in J(x)$. Since $x\in(V_{j_1})_{{3\over4}\delta}\bigcap\Ome$,
%$\lam_{j_1}(V_{j_1}\cap\Ome)$ is a subgraph and condition
%(\ref{ZZ}) is satisfied, it
%follows that $\{x-\eps\alpha_{j_1}\xi_{j_1},
%0<\alpha_{j_1}\le1\}\subset\Ome$. Next, let
%$j_2\in J(x),j_2\ne j_1.$ Since $x-\eps\alpha_{j_1}\xi_{j_1}\in
%(V_{j_2})_{{3\over4}\delta-\eps}\bigcap\Ome$ for all $\alpha_{j_1}\in(0,1]$,
%for similar reasons it follows that
%$$
%\{x-\alpha_{j_1}\xi_{j_1}-\alpha_{j_2}\xi_{j_2}, 0<\alpha_{j_1},
%\alpha_{j_2}\le1\}\subset\Ome
%$$
%and so on.
%
%Moreover, the same argument shows that
%\begin{equation}\label{A4}
%\{x-\sum\limits_{j\in J(x)}\alpha_j\xi_j,~0<\alpha_j<\infty\}
%\cap(\bigcap\limits_{j\in J(x)}V_j)\subset\Ome.
%\end{equation}

3. Assume that $0<\eps\le {\rm min}\,\{{A_3\over2},{\delta\over4}\}$
and $x\in\partial \Ome$. Then
\begin{equation}\label{Z}
d(T_\eps(x))\le|T_\eps x-x|=\eps\bigg|\sum\limits_{j=1}^s
\xi_j\psi_j(x)\bigg|\le\eps\sum\limits_{j=1}^s
\psi_j(x)=\eps.
\end{equation}
%Since $T_\eps(\Ome)\subset\Ome$
%$$d(T_\eps(\Ome),\partial\Ome)=d(\partial T_\eps(\Ome),\partial\Ome)
%=d(T_\eps(\partial\Ome),\partial\Ome)\le\eps,
%$$
%hence $\Ome\setminus T_\eps(\Ome)\subset\partial_\eps\Ome$ and
%$\Ome\setminus\partial_\eps\Ome\subset T_\eps(\Ome).$

Given $x\in\Ome$ there exists $\sigma>0$ such that $x\in
T_\eps(\Ome)$ for all $0\le\eps<\sigma$ and
$x\in\partial(T_\sigma(\Ome))= T_\sigma(\partial\Ome).$ Hence
$d(x)\le\sigma.$ So $d(x)>\eps$ implies $x\in T_\eps(\Ome)$ and
$\Ome\setminus\partial_\eps\Ome\subset T_\eps(\Ome).$

4. Let $d_j(x)$ be the distance of
$x\in\Ome\cap V_j$ from the boundary $\partial\Ome$ in the
direction of $\xi_j$, i.e.,
$d_j(x)=\phi_j(\overline{\lam_j(x)})-(\lam_j(x))_N.$
Then it is known that there exist $A_4>0$ such that for
all $j=1,...,s$ and all $x\in\Ome\cap V_j$
$$
A_4d_j(x)^{1\over\gamma}\le d(x)\le d_j(x).
$$

5. Let $x\in\Ome$ and $j\in J(x)$. Then, for $0<\eps\le{\delta\over4}$,
$T_\eps(x)\in\Ome\cap(V_j)_{\delta\over2}$ and
$$
d(T_\eps(x))^\gamma\ge A_4^\gamma d_j(T_\eps(x)).
$$
By Step 2
$$
[T_\eps(x),T_\eps(x)+\eps\psi_j(x)\xi_j]
\subset\{T_\eps(x)+\eps(\psi_j(x)-\alpha_j)\xi_j,0\le\alpha_j\le\psi_j(x)\}
$$$$
=\{x-\eps\alpha_j\xi_j-\eps\sum\limits_{i\in J(x),i\ne j}\xi_i\psi_i(x),
0\le\alpha_j\le\psi_j(x)\}\subset x-\overline{C(x,\eps)}\subset\Ome,
$$
hence
$$
d_j(T_\eps(x))\ge\eps\psi_j(x).
$$
Thus for all $j\in J(x)$
$$
d(T_\eps(x))^\gamma\ge A_4^\gamma\eps\psi_j(x)
$$
and
$$
sd(T_\eps(x))^\gamma\ge|J(x)|d(T_\eps(x))^\gamma
\ge A_4^\gamma\eps\sum\limits_{j\in J(x)}\psi_j(x)=A_4^\gamma\eps.
$$
Consequently
$$
d(T_\eps(x))\ge A_4s^{-{1\over\gamma}}\eps^{1\over\gamma},
$$
which implies that $T_\eps(\Ome)\subset\Ome
\setminus\partial_{A_1\eps^{1\over\gamma}}\Ome$
where $A_1=A_4s^{-{1\over\gamma}}.$

\begin{lemma}
Under the conditions of Lemma $18$ there exist $A_5,\eps_4>0$,
depending only on
$N,\gamma,M,$
$\delta,s,\{V_j\}_{j=1}^s,
\{\lambda_j\}_{j=1}^s,$
such
that for all $0<\eps\le \eps_4$
\begin{equation}\label{A2}
|\Ome\setminus T_\eps(\Ome)|\le A_5\eps.
\end{equation}
\end{lemma}

By (2)  $|\Ome\setminus \partial_\eps \Ome|\le A_6\eps^\gamma$
and $|\Ome\setminus \partial_{A_1\eps^{1\over\gamma}} \Ome|\le A_7\eps$.
Therefore the left inclusion of (22) immediately implies that
$|\Ome\setminus T_\eps(\Ome)|\le A_5\eps^\gamma$ but (\ref{A2}) makes a
stronger claim: estimate (\ref{A2}) has the same order in $\eps$ as the
estimate for $|\Ome\setminus \partial_{A_1\eps^{1\over\gamma}} \Ome|$.

In the proof of Lemma 19 the following property of regions satisfying the
cone condition will be used. We say that $C$ is a cone of size $\delta>0$ if
$$
C=\{(\bar x,x_D)\in \R^D:0<x_D<\delta,\bar x\in x_DK\}
$$
where $K$ is an open convex set in $\R^{D-1}$ such that
$0\in K\subset B(0,1).$

\begin{lemma}
Let $\delta>0$, $U,U'$ be bounded regions in $\R^{D-1}$ and
$U'\subset U_\delta$. Moreover, let $\phi,\psi:U\to\R$, $V=\{(\bar
x,x_D)\in\R^D:
\bar x\in U, \psi(\bar x)<x_D<\phi(\bar x)\}$ and
let $C\subset\R^D$ be a cone of size $\delta$. If
$$
(\bar x,\phi(\bar x))-C\subset V~~{\rm for ~all}~x\in U',
$$
then $\phi$ satisfies the
Lipschitz condition on $U'.$ Moreover, the Lipschitz constant depends
only on $K$.
\end{lemma}

The proof uses that $\{(\bar x,\phi(\bar x)+C\}\cap V=\O$ because otherwise
$(\bar x,\phi(\bar x))\in V.$

\underline{Proof of lemma 19}
1. For any subset $J\subset\{1,...,s\}$ put
$$
\tilde V_J=\{x\in\R^N:J(x)=J\},
$$
so $\R^N$ is the disjoint union of all $\tilde V_J:
~\R^N=\bigcap\limits_J\tilde V_J,$ and
\begin{equation}\label{*}
\Ome\setminus T_\eps(\Ome)=\bigcup\limits_J(\tilde
V_J\cap(\Ome\setminus T_\eps(\Ome))).
\end{equation}
        We will prove that
\begin{equation}\label{Q1}
\tilde V_J\cap(\Ome\setminus T_\eps(\Ome))\subset(\partial_{2\eps}\tilde V_J)
\cup \tilde V_J^{(\eps)},
\end{equation}
where
\begin{equation}\label{Q2}
\tilde V_J^{(\eps)}=\{\tilde V_J\cap\partial\Ome-A_J\eps\alpha\xi_J,
~0\le\alpha\le1\},
\end{equation}
$A_J$ is a certain positive number
and $\xi_J$ is a certain unit vector in $\R^N.$

2. First let $J=\{i\}.$ Then for all $x\in\tilde V_J$ we have $J(x)=\{i\}$ and
$\psi_i(x)=\sum\limits_{j=1}^s\psi_j(x)=1.$ Hence for all
$x\in\tilde V_J$
$$
T_\eps(x)=x-\eps\xi_i
$$
and for all $z\in T_\eps(\tilde V_J)$
$$
T^{-1}_\eps(z)=z+\eps\xi_i.
$$
Assume that $x\in\tilde V_J\cap(\Ome\setminus T_\eps(\Ome)).$
Since $x\not\in T_\eps(\Ome)$ and $T_\eps(x)\in T_\eps(\Ome)$, there is
a point $z\in\partial T_\eps(\Ome)=T_\eps(\partial\Ome)$ which lies in the
interval $(x-\eps\xi_i,x).$
Hence $y=T^{-1}_\eps(z)=z+\eps\xi_i\in\partial\Ome,$
$x\in(y,z)=(y,y-\eps\xi_i)\subset\tilde V_J^{(\eps)}$ and (\ref{Q1})
follows. (In this case the first entry of the union in the right-hand
side of (\ref{Q1}) can be omitted.)

3. Next let $J=\{i,j\}$. Then,
%8. Next we prove that there exist a vector $\xi_{ij}$ of length 1 and a
%positive number $A_{ij},$ depending only on $\xi_i$ and $\xi_j,$
%such that
%\begin{equation}\label{Q3}
%\tilde V_{ij}\cap(\Ome\setminus T_\eps(\Ome))
%\subset(\partial_{2\eps}\tilde V_{ij})
%\cup \tilde V_{ij}^{(\eps)}
%\end{equation}
%where
%\begin{equation}\label{Q4}
%\tilde V_{ij}^{(\eps)}
%=\{\tilde V_{ij}\cap\partial\Ome-A_{ij}\eps\alpha\xi_{ij},
%~0\le\alpha\le1\}.
%\end{equation}
for all $x\in\tilde V_J,$ $J(x)=\{i,j\},$
$$
T_\eps(x)=x-\eps(\psi_i(x)\xi_i+\psi_j(x)\xi_j).
$$
Let
$$
C_J=\{\alpha_1\xi_i+\alpha_2\xi_j,~\alpha_1,\alpha_2>0\},
$$
$$
C_J(\eps)=\{\eps(\alpha_1\xi_i+\alpha_2\xi_j),~\alpha_1,\alpha_2>0,~
\alpha_1+\alpha_2<1\}
$$
and let $x\in(V_i)_{3\delta\over4}\cap(V_j)_{3\delta\over4}\cap\Ome.$
Since $J\subset J(x)$, by Step 2 of the proof of Lemma 18 it follows that
\begin{equation}\label{B3}
V_i\cap V_j\cap(x-C_J)\subset V_i\cap V_j\cap(x-\overline{C(x)})\subset\Ome
\end{equation}
and
\begin{equation}\label{B2}
x-C_J(x)\subset x-\overline{C(x,\eps)}\subset V_i\cap V_j\cap\Ome.
\end{equation}

%and by Step 2 of the proof of Lemma 18
%\begin{equation}\label{B2}
%S_J(x)
%=\{x-\alpha_1\xi_i-\alpha_2\xi_j,0<\alpha_1,\alpha_2\le1\}\subset\Ome
%\end{equation}
%and
%\begin{equation}\label{B3}
%\{x-\alpha_1\xi_i-\alpha_2\xi_j,0<\alpha_1,\alpha_2<\infty\}\cap V_i
%\cap V_j\subset\Ome.
%\end{equation}
%(If $x\in\tilde V_J\cap\Ome$ then $\overline{S_{ij}(x)}\subset\Ome.)$

3.1. First assume that $\xi_i$ and $\xi_j$ are proportional.
Let $\xi_i=-\xi_j.$  If $x\in V_i\cap V_j\cap\partial\Ome,$ then
$x-\alpha\xi_i\in\Ome$ for small $\alpha>0$ and
$x-\alpha\xi_j\in\Ome$ for small $\alpha>0$. So
$x+\alpha\xi_i\in\Ome$ and $x-\alpha\xi_i\in\Ome$ for small $\alpha>0$
which contradicts condition 2) in the definition of a boundary of class
Lip $\gamma$.
Thus $\xi_i\ne-\xi_j,$ hence $\xi_i=\xi_j$ and
$$
T_\eps(x)=x-\eps(\psi_i(x)+\psi_j(x))\xi_i=x-\eps\xi_i.
$$
Similarly to Step 2 we obtain inclusion (\ref{Q1}) -- (\ref{Q2}) where
$\xi_J=\xi_i$ and $A_J=1.$

3.2.  Next assume that $\xi_i$ and $\xi_j$ are not proportional and set
$\xi_J=\frac{\xi_i+\xi_j}{|\xi_i+\xi_j|}.$
Let $\lam_J$ be a
rotation such that
$
\lam_J(\xi_J)=e_N
$ and the image of the plane spanned by $\xi_i$ and $\xi_j$ is the plane
spanned by $e_N,e_{N-1}$.
Inclusion (\ref{B3}) implies that
\begin{equation}\label{AB}
V_i\cap V_j\cap\{x-\alpha\xi_J,~\alpha>0\}
\subset\Ome.
\end{equation}
%hence
%$$
%\lam_{ij}(\Ome\cap\tilde V_{ij})=\{x_N<\phi_{ij}(\bar x),~\bar x\in G_{ij}\}
%\cap\lam_{ij}(\tilde V_{ij}),
%$$
%where $G_{ij}$ is the projection of $\lam_{ij}(\tilde V_{ij}\cap\partial\Ome)$
%onto the hyperplane $x_N=0.$ Moreover, inclusion (\ref{B3}) implies that the
%function $\phi_{ij}$ satisfies the Lipschitz condition in $x_{N-1}$ uniformly
%with respect to $x_1,...,x_{N-2}:$
%$$
%|\phi_{ij}(x_1,...,x_{N-2},x_{N-1})-\phi_{ij}(x_1,...,x_{N-2},y_{N-1})|
%\le L_{ij}|x_{N-1}-y_{N-1}|
%$$
%for $(x_1,...,x_{N-2},x_{N-1}),(x_1,...,x_{N-2},y_{N-1})\in G_{ij},$ where
%$L{ij}={\rm cot}\,{\alpha_{ij}\over2}$ and
%$\alpha_{ij}$ is the angle between the
%vectors $\xi_i$ and $\xi_j.$

3.3. Let
$$
H_J=\lam_J(V_i\cap V_j\cap\Ome),~~H'_J=\lam_J((V_i)_{3\delta\over4}
\cap (V_j)_{3\delta\over4}\cap\Ome)
$$
and $G_J,G_J'$ be the projections of $H_J,$ $H_J'$ respectively, onto
the hyperplane $x_N=0.$ Condition (\ref{AB}) implies that for all $\bar
x\in G_J$ there exist $\phi_J(\bar x),\psi_J(\bar x)$ such that $(\bar
x,\phi_J(\bar x))\in\Ome,~\psi_J(\bar x)<\phi_J(\bar x)$ and the intersection
of the line, parallel to $e_N$ and passing through $(\bar x,0)$, and $H_J$ is
$(\psi_J(\bar x),\phi_J(\bar x)).$ Hence
$$
H_J=\{(\bar x,x_N)\in\R^N:~\bar x\in G_J,~\psi_J(\bar x)<x_N<\phi_J(\bar
x)\}.
$$
Furthermore, by (\ref{B2}), for all $\bar x\in G_J'$
$$
(\bar x,\phi(\bar x))-\lam_J(C_J(\eps))\subset H_J.
$$
3.4. Next we apply, for fixed $x_1,...,x_{N-2},$ Lemma 20 where $D=2,$
$$
C=\lam_J(C_J(\eps))\subset\R^2_{e_{N-1},e_N}=\{(x_{N-1},x_N):
x_{N-1},x_N\in\R\},
$$
$\delta=\eps, K=\{x_{N-1}<\sin{\alpha_J\over2}\}$ ($\alpha_J$ is the
angle between $\xi_i$ and $\xi_j$), $U$ and $U'$ are the projections in
$\R^2_{e_{N-1},e_N}$ onto the line $x_N=0$ of the cross-sections of
$_J$, $G_J'$ respectively, in $\R^{N-1}_{e_1,...,e_{N-1}}$ by the line
parallel to $e_{N-1}$ and passing through $(x_1,...,x_{N-2},0)$,
$\phi=\phi_J, \psi=\psi_J$. It
follows that, for fixed $x_1,...,x_{N-2}$, the function $\phi_J$
satisfies the Lipschitz condition in $x_{N-1}$. Moreover, the Lipschitz
constant $L_J$ depends only on $K$, hence on $\alpha_J$, therefore is
independent of $x_1,...x_{N-2}$. (In fact $L_J=\cot{\alpha_J\over2}.$)
Thus, the
function $\phi_J$ satisfies the Lipschitz condition in $x_{N-1}$ uniformly
with respect to $x_1,...,x_{N-2}:$
$$
|\phi_J(x_1,...,x_{N-2},x_{N-1})-\phi_J(x_1,...,x_{N-2},y_{N-1})|
\le L_J|x_{N-1}-y_{N-1}|
$$
for all $x_1,...,x_{N-2},x_{N-1},y_{N-1}$ satisfying $(x_1,...,x_{N-2},
x_{N-1}),(x_1,...,x_{N-2},y_{N-1})\in G_J'.$

3.5. Assume that $x\in(\tilde V_J)_{2\eps}\cap(\Ome\setminus T_\eps(\Ome)).$
Since $x\not\in T_\eps(\Ome)$ and $T_\eps(x)\in T_\eps(\Ome)$, there is
a point $z\in\partial T_\eps(\Ome)=T_\eps(\partial\Ome)$ which lies in the
interval $[x,T_\eps(x)).$ Let $y=T_\eps^{-1}(z),$ then $y\in\partial\Ome.$
If $y\not\in \tilde V_J,$ then by (\ref{Z})
$$
d(x,\partial\tilde V_J)\le|x-y|\le|x-z|+|z-y|<|x-T_\eps(x)|
+ |y-T_\eps(y)|\le2\eps,
$$
which contradicts the assumption $x\in(\tilde V_J)_{2\eps}.$ So $y\in\tilde
V_J\cap\partial\Ome,$ hence $z\in y-C_J(\eps)$ where $C_J(\eps)$.

Let  $\lam_J(x)=\beta, \lam_J(y)=\eta, \lam_J(z)=\zeta$
and let $d_J(x)$ denote
the distance of $x\in\Ome$ from $\partial\Ome$ in the direction of the vector
$\xi_J,$
then
$$
\beta=(\beta_1,...,\beta_{N-2},\beta_{N-1},\beta_N),
\eta=(\beta_1,...,\beta_{N-2},\eta_{N-1},\eta_N),
$$$$
\eta_N=\phi_{ij}(\beta_1,...,\beta_{N-2},\eta_{N-1}),
\zeta=(\beta_1,...,\beta_{N-2},\zeta_{N-1},\zeta_N),
$$$$
d_J(x)\le\phi_J(\beta_1,...,\beta_{N-2},\beta_{N-1})-\beta_N.
$$
Since also $z\in[x,T_\eps(x))\subset x-\overline{C_J(\eps)},$  it follows that
$$
|\eta_{N-1}-\zeta_{N-1}|,|\eta_N-\zeta_N|\le|\eta-\zeta|=|z-T_\eps(z)|\le\eps,
$$$$
|\zeta_{N-1}-\beta_{N-1}|\le|z-x|\le|x-T_\eps(x)|\le\eps.
$$
Also $\zeta_N<\beta_N.$
Therefore
$$
d_J(x)\le\phi_J(\beta_1,...,\beta_{N-2},\beta_{N-1})
-\phi_J(\beta_1,...,\beta_{N-2},\eta_{N-1})
$$$$
+\phi_J(\beta_1,...,\beta_{N-2},\eta_{N-1})
-\beta_N
\le L_J|\beta_{N-1}-\eta_{N-1}|+\eta_N-\zeta_N
$$$$
\le L_J(|\beta_{N-1}-\zeta_{N-1}|+|\zeta_{N-1}-\eta_{N-1}|)+\eps
\le(2L_{ij}+1)\eps.
$$
Hence (\ref{Q1}) -- (\ref{Q2}) holds where $A_J=2L_{ij}+1.$

4. The argument for the cases in which
the number of elements in $J$ is greater than 2 is similar.
Let, for example, $J=\{i,j,k\}.$
If ${\rm dim~Span}\,\{\xi_i,\xi_j,\xi_k\}=1$, then $\xi_i=\xi_j=\xi_k$,
and we set $\xi_J=\xi_i$ and argue as in Step 3.1. If
${\rm dim~Span}\,\{\xi_i,\xi_j,\xi_k\}=2$, then we take any two linearly
independent vectors, say $\xi_i, \xi_j$, set $\xi_J=\frac{\xi_i+\xi_j}
{|\xi_i+\xi_j|}$ and argue as in Steps 3.2 -- 3.5. If ${\rm dim~Span}\,
\{\xi_i,\xi_j,\xi_k\}=3$, then we set $\xi_{ijk}=\frac{\xi_i+\xi_j+\xi_k}
{|\xi_i+\xi_j+\xi_k|}$ and argue as in Steps 3.2 -- 3.5. In this
case the appropriate function $\phi_J$ satisfies the Lipschitz
condition in $x_{N-1},x_{N-2}$ uniformly with respect to
$x_1,...,x_{N-3}.$ In Lemma 20 one can take  $K$ to be the
projection of $C_J(\eps)$ onto $x_N=0$ and
$\delta=\inf{\{x_N:(x_{N-3},x_{N-2},x_{N-1})\in C_J(\eps)\}}$.

5. By (\ref{*}) and (\ref{Q1})
$$
|\Ome\setminus T_\eps(\Ome)|\le\sum\limits_{J\subset\{1,...,s\}}
(|\partial_{2\eps}\tilde V_J|+
|\tilde V_J^{(\eps)}|).
$$
Since each summand does not exceed $\eps$ multiplied by a constant
depending only on $N,\gamma,M,\delta,s,\{V_j\}_{j=1}^s$ and
$\{\lam_j\}_{j=1}^s$, inequality (\ref{A2}) follows.
\vskip0.3cm
The proof of the following theorem can be adapted to any situation in
which  maps $\Lambda_\eps$ exist with properties (\ref{A5}) --(\ref{A7}).

\begin{theorem}
Let $N\ge2,0<\gamma\le1,M,\delta>0$ and integer $s\ge1$. Moreover, let
$\{V_j\}_{j=1}^s$ be a family of bounded open cuboids and
$\{\lambda_j\}_{j=1}^s$ a family of rotations. Suppose that
$\Ome_1\subset\R^N$ is a bounded region such that $\partial\Ome_1\in
{\rm Lip}\,(\gamma,M,\delta,s,\{V_j\}_{j=1}^s,$
$\{\lambda_j\}_{j=1}^s).$

Then for every integer $n\ge1$ there exist $b_{n,7}=b_{n,7}(\Ome_1),
\eps_{n,7}=\eps_{n,7}(\Ome_1)>0$
such that for all $0<\eps\le\eps_{n,7}$ and for all
bounded regions $\Ome_2$, for which $\partial\Ome_2\in {\rm
Lip}\, (\gamma,M,\delta,s,\{V_j\}_{j=1}^s,
\{\lambda_j\}_{j=1}^s)$ and
$\Ome_1\setminus\partial_\eps\Ome_1\subset\Ome_2 \subset\Ome_1$,
the inequality
\begin{equation}
(1-b_{n,7}\eps^\gamma)\lam_{n,1}
\le\lam_{n,2}\le(1+b_{n,7}\eps^\gamma)\lam_{n,1}
\end{equation}
holds.
\end{theorem}
\Proof
1. Since $\partial\Ome_1\in {\rm Lip}\,(\gamma,M,\delta,s,
\{V_j\}_{j=1}^s,
\{\lambda_j\}_{j=1}^s)$ it
follows that inequality (\ref{C0}) holds, where $q>2$ and $c_{13}>0$
depend only on $N,\gamma,M,\delta,s,\{V_j\}_{j=1}^s,
\{\lambda_j\}_{j=1}^s$ \cite{Glo},
\cite{Maz1}, \cite{Maz2}. Since also $\partial\Ome_2\in
{\rm Lip}\,(\gamma,M,\delta,s,
\{V_j\}_{j=1}^s,\{\lambda_j\}_{j=1}^s)$ inequality (\ref{C0}) holds with
$\Ome_2$ replacing $\Ome_1$ with the same $q$ and $c_{13}$.
Furthermore, by (2) there exist $A_8,\eps_5>0$, depending only on
$N,\gamma,M,\delta,s,\{V_j\}_{j=1}^s,\{\lambda_j\}_{j=1}^s$
such that for all $0<\eps\le\eps_5$
$$
|\partial_\eps\Ome_1|,|\partial_\eps\Ome_2|\le A_8\eps^\gamma.
$$

2. Let
$$
\Lambda_\eps=T_{({\eps\over A_1})^\gamma}~.
$$
Then by Lemmas 18 and 19 there exist $A_9,A_{10},A_{11},\eps_6>0$, depending only
on $N,\gamma,M,\delta,s,$
$\{V_j\}_{j=1}^s,\{\lambda_j\}_{j=1}^s)$, such that for $0<\eps\le\eps_6$
\begin{equation}\label{A5}
\Ome_3(\eps)=\Lambda_\eps(\Ome_1)\subset\Ome_1
\setminus\partial_\eps\Ome_1,
\end{equation}
\begin{equation}\label{A6}
|\Ome_1\setminus\Ome_3(\eps)|\le A_9\eps^\gamma
\end{equation}
and
\begin{equation}\label{A7}
\bigg|\frac{\partial\Lambda_{\eps i}}{\partial x_j}(x)-\delta_{ij}\bigg|
\le A_{10}\eps^\gamma,
\end{equation}
which implies that
\begin{equation}\label{A8}
{1\over2}\le 1-A_{11}\eps^\gamma\le
{\rm Jac}\,(\Lambda_\eps,x)
\le1+A_{11}\eps^\gamma,~~x\in\Ome_1.
\end{equation}
and the map
$\Lambda_\eps:\Ome_1\to \Lambda_\eps(\Ome_1)$ is one-to-one. Denote
$U_\eps=\Lambda_\eps^{-1}.$

3. Next we obtain a lower bound for $\mu_{n,3}$ by applying the
idea used in the proof of Lemma 13 based on the variational method.

If $L$ is an $(n+1)$-dimensional subspace of $L^2(\Ome_1)$, then
$\tilde L=\{f(U_\eps(x)),f\in L\}$ is an $(n+1)$-dimensional
subspace of $L^2(\Lambda_\eps(\Ome_1))$, and conversely, if $\tilde L$ is an
$(n+1)$-dimensional subspace of $L^2(\Lambda_\eps(\Ome_1))$, then
$L=\{g(\Lambda_\eps(x)),g\in \tilde L\}$ is an $(n+1)$-dimensional
subspace of $L^2(\Ome_1).$ Therefore
\begin{eqnarray*}
\mu_{n,3}
&=&\inf\limits_{\tilde L:\dim \tilde L=n+1}\sup\limits_{g\in \tilde
L}
\frac{\int_{\Ome_3(\eps)}|\nabla g|^2\rmd^Ny}
{\int_{\Ome_3(\eps)}|g|^2\rmd^Ny}\\ &=&\inf\limits_{L:\dim
L=n+1}\sup\limits_{f\in L}
\frac{\int_{\Lambda_\eps(\Ome_1)}|\nabla(f(U_\eps(y)))|^2\rmd^Ny}
{\int_{\Lambda_\eps(\Ome_1)}|f(U_\eps(y))|^2\rmd^Ny}.
\end{eqnarray*}
Note that
$$
|\nabla(f(U_\eps(y)))|^2=\sum\limits_{i=1}^N\bigg|\sum\limits_{k=1}^N
\bigg(\frac{\partial f}{\partial
x_k}\bigg)(U_\eps(y))\frac{\partial(U_\eps(y))_k}{\partial y_i}\bigg|^2
$$$$
=\sum\limits_{i=1}^N\bigg|\bigg(\frac{\partial f}
{\partial x_i}\bigg)(U_\eps(y))\bigg|^2
\bigg|\frac{\partial (U_\eps(y))_i}
{\partial y_i}\bigg|^2
$$$$
+\sum\limits_{i=1}^N{\sum\limits_{k,l=1}^N}\strut^i\bigg(\frac{\partial f}
{\partial x_k}\bigg)(U_\eps(y))\overline{\bigg(\frac{\partial f}
{\partial x_l}\bigg)(U_\eps(y))}
\frac{\partial (U_\eps(y))_k}
{\partial y_i}\overline{\frac{\partial (U_\eps(y))_l}
{\partial y_i}},
$$
where ${\sum}^i$ means that summation is taken with respect to such
$k,l$ that either $k\ne i$ or $l\ne i$.

Recall that
$$
\frac{\partial (U_\eps(y))_k}
{\partial
y_i}=\Delta_{ki}({\rm Jac}\,(\Lambda_\eps,U_\eps(y)))^{-1},
$$
where $(-1)^{k+i}\Delta_{ki}$ is the determinant obtained by
deleting $k$-th row and $i$-th column in the Jacobian
${\rm Jac}\,(\Lambda_\eps,U_\eps(y)))$. The Jacobi matrix of the map
$\Lambda_\eps$ has the form $I+(\frac\eps{A_1})^\gamma B$,
where $I$ and $B$ are defined
in the proof of Lemma 18. Hence $(-1)^{k+i}\Delta_{ki}$ is the
determinant of the matrix $I_{ki}+(\frac\eps{A_1})^\gamma B_{ki},$
where  $I_{ki}$
and $B_{ki}$ are obtained by deleting $k$-th rows and $i$-th
columns in matrices $I$, $B$ respectively. Since $|I_{ii}|=1$
and $|I_{ki}|=0$ if $k\ne i$ and the elements of the matrix $B$
are uiformly bounded, there exists $A_{12},\eps_7>0,$ depending only on
$N,\delta$ and $s$, such that for all $0<\eps\le\eps_7$
and $y\in \Lambda_\eps(\Ome_1)$
$$
{1\over2}\le 1-A_{12}\eps^\gamma \le\bigg|\frac{\partial (U_\eps(y))_i}
{\partial y_i}\bigg|\le1+A_{12}\eps^\gamma, ~~i=1,...,N,
$$
and
$$
\bigg|\frac{\partial (U_\eps(y))_k}
{\partial y_i}\bigg|\le A_{12}\eps^\gamma ,~~i,k=1,...,N,~~k\ne i.
$$
Consequently there exists $A_{13}>0,$ depending only on
$N,\delta$ and $s$, such that for all $0<\eps\le\eps_7$
and $y\in \Lambda_\eps(\Ome_1)$
$$
(1-A_{13}\eps^\gamma)|(\nabla f)(U_\eps(y))|^2
\le |\nabla (f(U_\eps(y)))|^2
\le(1+A_{13}\eps^\gamma)|(\nabla f)(U_\eps(y))|^2
$$
Therefore, by changing variables: $y=\Lambda_\eps(x)$ and taking into
account inequality (\ref{A7}), we have
$$
\mu_{n,3}\ge(1-A_{13}\eps^\gamma)\inf\limits_{L:\dim L=n+1}\sup\limits_{f\in L}
\frac{\int_{\Lambda_\eps(\Ome_1)}|(\nabla f)(U_\eps(y))|^2\rmd^Ny}
{\int_{\Lambda_\eps(\Ome_1)}|f(U_\eps(y))|^2\rmd^Ny}
$$$$
=(1-A_{13}\eps^\gamma)\inf\limits_{L:\dim L=n+1}\sup\limits_{f\in L}
\frac{\int_{\Ome_1}|(\nabla f)(x)|{\rm Jac}\,(\Lambda_\eps,x)|\rmd^Nx}
{\int_{\Ome_1}|f(x)|^2|{\rm Jac}\,(\Lambda_\eps,x)|\rmd^Nx}
$$$$
\ge(1-A_{13}\eps^\gamma)(1-A_{11}\eps^\gamma)(1+A_{11}\eps^\gamma)^{-1}
\inf\limits_{L:\dim L=n+1}\sup\limits_{f\in L}
\frac{\int_{\Ome_1}|\nabla f|^2\rmd^Nx}
{\int_{\Ome_1}|f|^2\rmd^Nx}.
$$
Hence, finally, there exist $b_9,
\eps_9>0$,
depending only on
$N,\gamma,M,$
$\delta,s,\{V_j\}_{j=1}^s,
\{\lambda_j\}_{j=1}^s,$
 such that for $0<\eps\le\eps_9$
$$
\mu_{n,3}\ge(1-\eps^\gamma b_9)\lam_{n,1}.
$$

Now the theorem follows by taking into account (\ref{A5}) and (\ref{A6}), and
applying Corollary 16.

\vskip 2cm
%%%%%%%%%%%%%%%%%%%%%%%%%
{\bf Acknowledgments } We would like to thank R Banuelos, W D Evans,
Y Safarov and M van den Berg for useful comments.
\newpage
%%%%%%%%%%%%%%%%%%%%%%%

\par
%%%%%%%%%%%%%%%%%%%%%%%%%%%%%%%%%%%%%%%%
\vskip 0.5in
\begin{tabbing}
School of Mathematics \qquad\qquad\qquad\qquad\= Department of Mathematics\\
Cardiff University \>Kings College's\\
23 Senghennydd Road\> Strand\\
Cardiff CF24 4YH \> London WC2R 2LS\\
Wales\>England\\
\>\\
e-mail: Burenkov@cardiff.ac.uk\>
e-mail: E.Brian.Davies@kcl.ac.uk\\
\end{tabbing}

\end{document}